\documentclass[11pt,twoside]{article}
\usepackage{amsmath,amsfonts,latexsym,amscd}

\scrollmode

\usepackage{amsthm}

\swapnumbers
\theoremstyle{plain}
\newtheorem{theorem}{Theorem}[section]
\newtheorem{lemma}[theorem]{Lemma}
\newtheorem{corollary}[theorem]{Corollary}
\newtheorem{proposition}[theorem]{Proposition}

\theoremstyle{definition}
\newtheorem{definition}[theorem]{Definition}

\newtheorem{situation}[theorem]{Situation}

\theoremstyle{remark}
\newtheorem{remark}[theorem]{Remark}

\usepackage{amsmath,amsfonts,latexsym}

\newcommand{\R}{\mathbb{R}}
\newcommand{\reals}{\mathbb{R}}
\newcommand{\C}{\mathbb{C}}

\newcommand{\N}{\mathbb{N}}

\newcommand{\1}{\mathbf{1}}  
\DeclareMathOperator{\id}{id}
\newcommand{\boundary}[1]{\partial#1}
\newcommand{\boundedops}{\mathcal{B}}
\newcommand{\bundleover}{\!\!\downarrow\!\!}
\newcommand{\abs}[1]{\left\lvert#1\right\rvert} 
\newcommand{\norm}[1]{\left\lVert#1\right\rVert}

\newcommand{\tensor}{\otimes}

\DeclareMathOperator{\supp}{supp}   
\DeclareMathOperator{\im}{im}      
\DeclareMathOperator{\Hom}{Hom}    

\DeclareMathOperator{\tr}{tr}
\DeclareMathOperator{\res}{res}  
\DeclareMathOperator{\fl}{f\/l}  
\DeclareMathOperator{\coker}{coker}

\DeclareMathOperator{\ind}{ind}

\newcommand{\forget}[1]{}

{\catcode`@=11\global\let\c@equation=\c@theorem}

\allowdisplaybreaks[2]

\newcommand{\citediss}{Schick(1996)url}

\newcommand{\bvp}[1]{\mathcal{#1}} 
\newcommand{\Neumann}[1]{\mathcal{#1}}
\newcommand{\BdM}[1]{\mathcal{#1}} 
\renewcommand{\Sp}{\rm{Sp}}

\begin{document}

\title{$L^2$-index theorem for boundary manifolds}
\author{Thomas Schick} 
     
\date{}        
\maketitle

\begin{abstract}
Suppose $M$ is a compact manifold with boundary $\boundary M$. Let
$\tilde M\bundleover M$ be a normal covering with covering group
$\Gamma$. Suppose $(A,T)$ is an elliptic differential boundary value
problem on $M$ with lift $(\tilde A,\tilde T)$ to $\tilde M$. Then the
von Neumann dimension $\dim_\Gamma$ of kernel and cokernel of this
lift are defined. The main result of this
paper is: these numbers are finite, and their
difference, by definition the von Neumann index $\ind_\Gamma(\tilde
A,\tilde T)$, equals the index of $(A,T)$.  In this way, we extend the
classical $L^2$-index
theorem of Atiyah to manifolds with boundary.

\smallskip
MSC-classification: 58G12 (primary); 58G03 (secondary)
\end{abstract}

\section{Introduction}
In this paper, we study elliptic differential boundary value problems
on coverings of compact manifolds. Let $M$ be a compact Riemannian
manifold with boundary $\boundary M$. Suppose $E,F\bundleover M$ and
$Y\bundleover \boundary M$ are Riemannian vector bundles. Let
$A:C^\infty(E)\to C^\infty(F)$ be a differential operator and
$T:C^\infty(E)\to C^\infty(Y)$ a differential boundary operator so
that the pair $\bvp{P}:=(A,T)$ is elliptic.  Define 
\[\begin{split} \ker\bvp{P}&:=\{f\in L^2(E);\; f\in C^\infty, Af=0=Tf\}\qquad\text{ and}\\
  \coker\bvp{P}&:=\{ (F,f)\in L^2(F)\oplus L^2(Y);\\
   &\qquad (F,A\varphi)_{L^2(F)}+(f,T\varphi)_{L^2(Y)}=0\;\;\forall\varphi\in C_0^\infty(E)\}.
\end{split} \]
The classical theory of elliptic boundary value problems states that
the dimensions of kernel and cokernel are finite and studies
$\ind(\bvp{P}):=\dim\ker\bvp{P} -\dim\coker\bvp{P}$. The index theorem
(recalled below) provides deep connections between topological, geometrical and analytical properties of the manifold.

Suppose $\tilde M\bundleover M$ is a normal covering of $M$ with
deck transformation group $\Gamma$. Pull the bundles back to $\tilde M$ and
lift the operators and metrics. We use the convention that
corresponding objects on $\tilde M$ have the same notation decorated
with an additional tilde. Note that $\Gamma$ operates on the bundles,
their sections and that $\tilde\bvp{P}=(\tilde A,\tilde T)$ is
$\Gamma$-equivariant. Define the kernel and cokernel of
$\tilde\bvp{P}$ literally in the same way as for $\bvp{P}$. They are in
general infinite dimensional. But $\ker(\tilde\bvp{P})$ and
$\coker(\tilde\bvp{P})$ are Hilbert modules over the group von Neumann
algebra $\Neumann{N}(\Gamma)$. (For von Neumann
algebras and Hilbert modules compare
\cite{Dixmier(1969),Lott-Lueck(1995)}.) For these Hilbert modules, a
normalized dimension $\dim_\Gamma$ with values in $[0,\infty]$ is
defined. It vanishes exactly if the module is trivial, it is additive
under direct sums, and 
\begin{equation}
  \label{finitemult}
  \abs{\Gamma}<\infty\quad\implies
\dim_\Gamma=\frac{1}{\abs{\Gamma}}\dim_{\C} .
\end{equation}

The following is the main result of this paper:
\begin{theorem}\label{ind_theorem} 
In the situation described above, 
$\dim_\Gamma\ker(\tilde\bvp{P}) <\infty$, $\dim_\Gamma\coker(\tilde
\bvp{P}) <\infty$ and
\[
\ind_\Gamma(\tilde\bvp{P}):=\dim_\Gamma\ker(\tilde\bvp{P})-\dim_\Gamma\coker(\tilde\bvp{P})
= \ind(\bvp{P}) . \]
Remarkably, $\ind_\Gamma(\tilde\bvp{P})$,  the difference of
two reals, is an integer.
\end{theorem}

The theorem is particularly interesting because for $\ind(\bvp{P})$ on
$M$ a well known {\em purely topological} expression exists, compare
Atiyah/Bott 
\cite[Theorem 2]{Atiyah-Bott(1964)} or Atiyah's \cite[Appendix
I]{Palais(1965)}: every elliptic boundary value problem $(A,p)$
defines a $K$-theoretic symbol class $[\sigma(A,p)]\in K(B(M),S(M)\cup
B(M)|_{\boundary M})$ (here $B(M)$ and $S(M)$ are the disc bundle and
sphere bundle of $TM$). As usual  one assigns
to this symbol class the topological index, and it coincides with the
analytical index. It can be computed cohomologically as (compare
\cite{Atiyah-Bott(1964)})
\[ \ind_t(A,p)= \int_{S(M)}ch(A)\pi^*\mathcal{T}(M) +
\int_{B(M)|_{\boundary M}} ch(A,p)\pi^*\mathcal{T}(M), \]
where $\pi:TM\to M$ is the projection, $\mathcal{T}(M)$ is a Todd class of $M$ and $ch$ the Chern character.

\begin{corollary}[of Theorem \ref{ind_theorem}]\label{mult}
The index of elliptic differential boundary value problems is multiplicative under finite coverings.
\end{corollary}
\begin{proof}
This follows from the multiplicativity \eqref{finitemult} of $\dim_\Gamma$.
\end{proof}

In Theorem \ref{ind_theorem} we can replace $\coker(\tilde\bvp{P})$ with
the kernel of an adjoint boundary value problem  by Theorem \ref{cokernel}.
Sometimes it is  easier to deal with
kernels. As an application we
compute the Euler characteristic of $M$ in terms of $L^2$-harmonic forms on
$\tilde M$ in Theorem \ref{Euler}. Dodziuk \cite{Dodziuk(1979)} and
Donnelly/Xavier \cite{Donnelly-Xavier(1984)} have computed the sign of
the Euler characteristic of closed negatively curved manifolds
in this way. An extension to manifolds with boundary is given in
\cite[Section 6]{\citediss}.

\smallskip
Our index theorem is the generalization of Atiyah's $L^2$-index
theorem \cite{Atiyah(1976)} to manifolds with boundary. The proof is
along the lines of Atiyah's proof. In order to deal with boundary
value problems, we  replace the calculus of pseudo-differential
operators by the Boutet de Monvel  calculus. We also try to clarify the
exposition. For this reason, in Section \ref{sec_trace} we introduce
and study thoroughly traces for endomorphisms of arbitrary Hilbert
modules. We use
the theory of Sobolev spaces to simplify the work with regularizing
operators and especially with their traces. An important result, which
should be valuable also in other contexts, is:
\begin{theorem}\label{L2_Rellich0} (compare Theorem \ref{L2_Rellich})\\
If $r>\dim M/2$, the inclusion of Sobolev spaces
$ H^{s+r}(\tilde M) \hookrightarrow  H^s(\tilde M)$
is a $\Gamma$-trace class operator. 
\end{theorem}

The  idea for the proof of the index theorem is: to $\bvp{P}$
construct an inverse $\bvp{Q}$ (modulo smoothing operators) in the BdM
calculus which  can be lifted to $\tilde M$, i.e.\
$\bvp{P}\bvp{Q}=\1-\bvp{S}_1$, $\bvp{Q}\bvp{P}=\1-\bvp{S}_0$ and
$\tilde\bvp{P}\tilde\bvp{Q}=\1-\tilde\bvp{S}_1$,
$\tilde\bvp{Q}\tilde\bvp{P}=\1-\tilde\bvp{S}_0$. Then the following
two results prove the theorem:
\begin{itemize}
\item $\ind_\Gamma(\tilde\bvp{P})=\Sp_\Gamma\tilde\bvp S_0-\Sp_\Gamma\tilde\bvp S_1$ (and the corresponding formula on the base with $\Gamma=\{1\}$)
\item For lifts of smoothing operators, we have $\Sp_\Gamma\tilde\bvp{S}=\Sp\bvp{S}$.
\end{itemize}

Note that our index theorem does not generalize the
Atiyah-Patodi-Singer index theorem
\cite{Atiyah-Patodi-Singer(1975a)}. They deal with a
specific non-local boundary condition for Dirac type operators. There is also an
$L^2$-version of this type of index theorem, proved by Ramachandran
\cite{Ramachandran(1993)}. He deals with Dirac type operators
and the APS-boundary conditions. Contrariwise, our result is valid for
arbitrary elliptic {\em differential} boundary value problems, but we only
deal with local boundary conditions. In particular, we can not handle
the signature.

This work is part of the Dissertation \cite{\citediss} of the
author. I thank my advisor Prof.~Wolfgang L\"uck for his
constant support. 

Throughout the paper, we use the following notation:
\begin{definition}
For $c>0$ we define
\[ a\stackrel{c}{\le} b\quad\iff\quad a\le c\cdot b \]
and similarly $a\stackrel{c}{<}b$, \ldots.
In a longer chain of inequalities, the same symbol (e.g. $c$)  may be used for different constants.

If not stated otherwise, $H$ is a Hilbert space, $\boundedops(H_1,H_2)$ denotes the bounded operators from $H_1$ to $H_2$, $\boundedops(H)=\boundedops(H,H)$.
$M$ is a compact smooth manifold of dimension $m$ with boundary
$\boundary M$ and $E,F\bundleover M$, $Y\bundleover \boundary M$ are
vector bundles.
\end{definition} 


\section{Traces for $\Neumann{N}(\Gamma)$-module morphisms}\label{sec_trace}

In this section, we introduce the von Neumann trace $\tr_\Gamma$ for
morphisms between (not necessarily finite) Hilbert modules
over  the von Neumann algebra $\Neumann{N}(\Gamma)$ of a discrete
group $\Gamma$. To gain greater
flexibility, we introduce the concept of $\Gamma$-trace class
operators also if domain and range are different.

First note that on $\Neumann{N}(\Gamma):=\boundedops(l^2\Gamma)^\Gamma$  we have the canonical finite trace $\tr_\Gamma(a)=(a(e),e)$. Moreover, on a Hilbert space $H$ the trace $\Sp(A)=\sum_i(Ah_i,h_i)$ exists (($h_i$) an orthonormal basis of $H$). 
\begin{definition}
This yields a $\Gamma$-trace, called $\Sp_\Gamma$, on the $\Gamma$-operators on $l^2(\Gamma)\tensor H$ which is defined by
\[ \Sp_\Gamma(a\tensor A)=\tr_\Gamma(a)\cdot \Sp(A) .\]
Note that this makes sense only for positive operators and for the operators in the $\Gamma$-trace class ideal (defined as usual, compare \cite[chapter I]{Dixmier(1969)}). We also have the $\Gamma$-Hilbert Schmidt (HS) operators defined by
\[ f\in\boundedops(l^2\Gamma\tensor H)^\Gamma\text{ is $\Gamma$-HS}\quad\iff\quad \Sp_\Gamma(f^*f)<\infty .\]
\end{definition}

Now we handle the general case:
Remember that a {\em Hilbert $\Neumann{N}(\Gamma)$-module} is a Hilbert space $V$ with  left $\Gamma$-action so that an isometric embedding $V\hookrightarrow l^2(\Gamma)\tensor H$ exists which is compatible with the $\Gamma$-actions.
If $V,W$ are two Hilbert $\Neumann{N}(\Gamma)$-modules, a bounded linear map $f:V\to W$ which is compatible with the $\Gamma$-action is called an {\em$\Neumann{N}(\Gamma)$-module morphism}.

\begin{definition}\label{def_tr}
Let $V_k$ be Hilbert $\Neumann{N}(\Gamma)$-modules with isometric $\Gamma$-em\-bed\-dings $i_k:V_k\hookrightarrow l^2(\Gamma)\tensor H_k$. Let $p_k:=i_k^*$ be the corresponding projections ($k=1,2$). Let $f:V_1\to V_2$ be a Hilbert $\Neumann{N}(\Gamma)$-module morphism.

We call $f:V_1\to V_2$ {\em $\Gamma$-Hilbert Schmidt ($\Gamma$-HS)} if $\Sp_\Gamma(i_1 f^*fp_1)<\infty$,
and we denote it {\em $\Gamma$-trace class ($\Gamma$-tr)} if
$\Gamma$-HS morphisms $f_1:V_1\to V_3$ and
$f_2:V_3\to V_2$ exist so that $f=f_2f_1$ .

If $V_1=V_2$ and $f$ is $\Gamma$-tr we set
\[ \Sp_\Gamma(f):=\Sp_\Gamma(i_1 f p_1) .\]
\end{definition}

The following basic properties show that this is well defined.

\begin{theorem}\label{tr_prop}
Let $f: V_1\to V_2$, $g: V_2\to V_3$, $e: V_0\to V_1$ be Hilbert $\Neumann{N}(\Gamma)$-module morphisms. Then
\begin{enumerate}
\item\label{adjoint}\label{tr_descr} $f$ $\Gamma$-tr $\iff$ $f^*$
  $\Gamma$-tr  $\iff$ $\abs{f}$ $\Gamma$-tr; $\quad$  $f$ $\Gamma$-HS $\iff$ $f^*$ $\Gamma$-HS
\item\label{HS_comp} $f$ $\Gamma$-HS $\implies$ $gf$, $fe$ $\Gamma$-HS
\item\label{tr_comp} $f$ $\Gamma$-tr $\implies$ $gf$, $fe$ $\Gamma$-tr
\item\label{C_weak} $f$ $\Gamma$-tr and $V_1=V_3$ $\implies$ $g\mapsto\Sp_\Gamma(gf)$ is ultra-weakly continuous.
\item\label{tr_b}\label{tr_HS} $V_1=V_3$ and either $f$ $\Gamma$-tr or $f, g$
  $\Gamma$-HS $\implies \Sp_\Gamma(gf)=\Sp_\Gamma(fg)$
\item\label{HS_prod}\label{tr_prod} If $V_{1,2}=l^2(\Gamma)\tensor H$,
  $a$ is $\Gamma$-HS  and $B\in\boundedops(H)$ is HS, then $f=a\tensor
  B$ is $\Gamma$-HS.
If $a$ is $\Gamma$-tr and $B$ is trace class, then $f$ is $\Gamma$-tr with $\Sp_\Gamma(f)=\tr_\Gamma(a)\Sp(B)$. 
\end{enumerate}
\end{theorem}
\begin{proof} This is a consequence of the corresponding properties of
  $\tr_\Gamma$ and $\Sp$, compare \cite[9.13]{\citediss}.
\end{proof}

As usual, armed with a $\Gamma$-trace we define the $\Gamma$-dimension:
\begin{definition}
Let $V$ be a  Hilbert $\Neumann{N}(\Gamma)$-module. Then
\[ \dim_\Gamma(V):= \Sp_\Gamma(\id_V) \in [0,\infty]. \]
\end{definition}

We now come to an important result, which is essentially proved in
Atiyah's paper \cite{Atiyah(1976)}. He does not state it in full
generality, but his proof works nearly literally, and can also be found
in \cite[9.16]{\citediss}.
\begin{proposition}\label{equaltr}
Suppose $V,W$ are Hilbert $\Neumann{N}(\Gamma)$-modules. Let $T_0:V\to V$ and $T_1:W\to W$ be bounded $\Gamma$-morphisms which are $\Gamma$-tr. Let $D:V\to W$ be a closed operator with domain $\mathcal{D}(D)$ which commutes with the action of $\Gamma$. Especially, we require that $\mathcal{D}(D)$ is $\Gamma$-invariant and dense. Suppose
\[ DT_0=T_1D;\quad \ker D\subset\ker T_0;\quad \ker D^*\subset\ker T_1^*.\]
\[ \implies \Sp_\Gamma(T_0) = \Sp_\Gamma(T_1) \]
\end{proposition}


\section{$L^2$-Rellich lemma}\label{sec_Rellich}

Let $M$ be a compact $m$-dimensional manifold with boundary $\boundary M$ (possibly empty). Let $\tilde M$ be a normal covering of $M$ with covering group $\Gamma$ (acting by isometries). Let $E\bundleover M$ be a vector bundle with pullback $\tilde E\bundleover \tilde M$.

There is a natural way to define Sobolev spaces on $\tilde M$:
\begin{definition}
Choose a finite covering of $M$ by charts $\kappa_i$ with subordinate partition of unity $\varphi_i$ so that $E$ is trivial over the domain of $\kappa_i$ with trivialization $t_i$. Lift charts, partition of unity and trivializations to $\tilde M$. Then we define the Sobolev norm $\abs{\cdot}_{H^s}$ by
\[ \abs{\sigma}_{H^s}:= \sum{\gamma\in\Gamma}\sum_{i}\abs{\tilde t_i
\circ(\tilde\varphi_i\cdot\gamma^*\sigma)\circ\tilde
\kappa_i^{-1}}_{H^s(\R^m)} \qquad \sigma\in C^\infty_0(\tilde E). \]
The Sobolev space $H^s(\tilde E)$ is defined as the completion of
$C^\infty_0(\tilde E)$ with respect to this norm. The inner product
does depend on the choices, but not the topology.
\end{definition}

We will show in this section that $H^s(\tilde E)$ is a Hilbert $\Neumann{N}(\Gamma)$-module and that the inclusion $H^{s+r}(\tilde E)\hookrightarrow H^s(\tilde E)$ is $\Gamma$-HS for $r>m/2$. 

Let $W$ be the double of $M$ with reflection $\fl: W\to W$. Let $X\bundleover W$ be the double of $E$. The reflection $\fl$ extends as a bundle map to $X$. Construct similarly $\tilde W$ and $\tilde X$. Then $\tilde W$ is a normal covering of $W$ with covering group $\Gamma$. Again we denote the reflection $\fl$. 

\begin{lemma}\label{ext}
Fix $s\in\R$. There exists a bounded $\Gamma$-equivariant extension
map $e: H^s(\tilde M)\to H^s(\tilde W)$, i.e.\ $e(f)|_{\tilde
  M}=f\;\;\forall f\in H^s(\tilde M)$. The restricition map is also
$\Gamma$-equivariant and bounded.\\
The corresponding statement holds for $\tilde E$.
\end{lemma}
\begin{proof}
The proof can be found in \cite[p.~27]{Schick(1994)}. One uses a $\Gamma$-invariant covering of $\tilde M$ by charts and  the corresponding extension map on Euclidian space (Taylor \cite[I.5.1]{Taylor(1981)}).
\end{proof}

Suppose $U\subset\tilde M\subset \tilde W$ is a fundamental domain for
the covering $p:\tilde M\to M$. This means that $U$ is open, $p|_U$ is
injective and $M-p(U)$ is a set of measure zero. Choose $U$ so that
its closure is compact, and choose a compact submanifold with boundary
$T\subset \tilde W$ of codimension zero, so that $U\cup\fl(U)\subset
T$ and so that the interior of $T$ is mapped surjectively onto $W$.

\begin{lemma}
Suppose $s\in\R$.
The map $p$ defined by the composition
\[\begin{CD} H^s(\tilde M)@>e>> H^s(\tilde W)@>{\bar p}>> l^2(\Gamma)\tensor H^s(T)\\
\in && \in && \in\\
                f @>>> ef @>>> \sum_{g\in\Gamma} g\tensor \underbrace{g^{-1}(ef)}_{=g^*(ef)}|_T
  \end{CD}\]
is $\Gamma$-equivariant, and there exist $C_{1,2}>0$ so that
\[ \abs{f}_{H^s(\tilde M)} \stackrel{C_1}{\le} \abs{pf}_{l^2(\Gamma)\tensor H^s(T)}\stackrel{C_2}{\le} \abs{f}_{H^s(\tilde M)}. \] 
In particular, $H^s(\tilde M)$ (with the pull back norm under $p$) is
a Hilbert $\Neumann N(\Gamma)$-module. 
The corresponding statement holds for $\tilde E$.
\end{lemma}
\begin{proof}
By Lemma \ref{ext}, $e$ has the required properties. It remains to consider $\bar p$. Obviously, $\bar p$ is $\Gamma$-equivariant.

Because $\Gamma$ is discrete and $T$ is compact, it meets only finitely many, say $N$,  of its translates $\{gT\}_{g\in\Gamma}$.

By definition, $\abs{\sum g\tensor f_g}^2_{l^2(\Gamma)\tensor H^s(T)}=\sum\abs{f_g}^2_{H^s(T)}$.
To show that $\bar p$ is bounded
let $\{U_i\}_{i=1,\dots N}$ be open subsets of $\tilde W$ which cover  $T$ so that the
covering projection maps each $U_i$ injectively to $W$. Choose submanifold charts $\kappa_i$ for $(U_i,U_i\cap T)$ and functions $0\le \varphi_i\le 1$ with compact support in $U_i$ so that $\sum_i\varphi_i=1$ on $T$.
Recognize that for every single $i$ we can extend $(U_i,\varphi_i,\kappa_i)$ to a corresponding collection $(U^i_{\alpha,\gamma},\varphi^i_{\alpha,\gamma},\kappa^i_{\alpha,\gamma})_{\alpha,\gamma}$ which can be used to compute Sobolev norms on $\tilde W$.
The norm will depend on the data (hence on $i$), but all the norms are equivalent. Therefore for $f\in H^s(\tilde W)$

\begin{align*}
\abs{\bar p f}_{l^2(\Gamma)\tensor H^s(T)}^2 &= \sum_{i=0}^N\sum_{\gamma\in\Gamma} 
\abs{\varphi_i\gamma^*f\circ\kappa_i^{-1}}_{H^s(\R^m_{\ge 0})}^2 \\
&\le \sum_i\sum_\gamma\sum_{\alpha=1}^{N_i}\abs{(\varphi_{\alpha,\gamma}^i f)\circ(\kappa^i_{\alpha,\gamma})^{-1}}_{H^s(\R^m)}^2\\
\intertext{(since we have more and larger summands)}
&\stackrel{NC}{\le} \abs{f}^2_{H^s(\tilde W)} .
\end{align*}
On the other hand (fix $i$)
\begin{align*}
 \abs{f}^2_{H^s(\tilde W)} &= \sum_{\alpha=1}^{N_i}\sum_\gamma
\abs{(\varphi^i_{\alpha,\gamma}f)\circ(\kappa^i_{\alpha,\gamma})^{-1}
}^2_{H^s(\R^m)}\\
\intertext{(choose $U^i_{\alpha,\gamma}$ 
so small that each of them lies in the interior of some translate of
$T$. Then we
can for every fixed $\alpha$ add more positive summands to get (up to norm
equivalence) $\abs{\cdot}_{l^2(\Gamma)\tensor H^s(T)}$. Therefore:)}
&\stackrel{CN_i}{\le} \abs{f}_{l^2(\Gamma)\tensor H^s(T)}  .
\end{align*}
The computations for $\tilde E$ are similar, but notationally more complicated.
\end{proof}

\begin{theorem}\label{L2_Rellich}
Suppose $s,r\in\reals$. The inclusion $\tilde{i}: H^{s+r}(\tilde
E)\to H^s(\tilde E)$ is $\Gamma$-HS if $r>m/2$, and is $\Gamma$-tr if $r>m$.
\end{theorem}
\begin{proof}
Let $X\bundleover W$ be the double of $E$. The following diagram commutes by the geometric definition of $p$:
\[\begin{CD}
 H^{s+r}(\tilde E) @>{p_{s+r}}>> l^2(\Gamma)\tensor H^{s+r}(\tilde X|_T)\\
 @V{\tilde{i}}VV    @VV{\1\tensor i}V\\
 H^s(\tilde E) @>{p_s}>> l^2(\Gamma)\tensor H^s(\tilde X|_T) .
\end{CD}\]
Remember that we have equipped $H^s(\tilde E)$ with the Hilbert space structure which makes $p$ an isometric embedding, therefore $p^*p=\1$. This yields
\[ \tilde{i} = p^*_s p_s\tilde{i}=p^*_s (\1\tensor i) p_{s+r} \]
Now we apply Properties (\ref{HS_comp}) and (\ref{HS_prod}) of Theorem
\ref{tr_prop}, together with the classical result that for bundles
over compact manifolds the inclusion
$H^{s+r}\hookrightarrow H^s$ is HS if $r> m/2$.
The second conclusion is an immediate corollary of the first.
\end{proof}


\section{Boutet de Monvel calculus}\label{sec_BdM}

The {\em Boutet de Monvel (BdM)} calculus is a tool to deal with
boundary value problems. It generalizes the calculus of
pseudo-differential operators on manifolds without boundary. We will
not go into the 
details, but only cite the results which are essential for our
applications.

The main point of the Boutet de Monvel calculus is the introduction of
an algebra of operators which includes the boundary value problems we
want to study and also their inverses. The first observation we have
to make is that we naturally have to consider matrices of operators.

For us the following is important: every elliptic boundary value
problem has a parametrix (an inverse modulo smoothing operators). And
every BdM operator is up to a smoothing operator nearly local.

\begin{definition}\label{def_BdM}
Let $M$ be a manifold with boundary $\boundary M$. Let $E,F\bundleover M$ be vector bundles over $M$, $X,Y\bundleover \boundary M$ bundles over the boundary. A BdM operator $\BdM{P}$ has the shape
\[ \BdM{P} = \begin{pmatrix} A+G & K\\
             T & p\end{pmatrix}      
             : \begin{matrix}C_0^\infty(E)\\ \oplus
\\C_0^\infty(X)\end{matrix} \to \begin{matrix}C^\infty(F)\\ \oplus\\
C^\infty(Y)\end{matrix} ,
             \]
where $A$ is a pseudo-differential operator \em{(pdo)} with the transmission property on $M$, $p$ is a pdo on $\boundary M$. $T: C^\infty(E)\to C^\infty(Y)$ is a trace operator, $K:C_0^\infty(X)\to C^\infty(F)$ a potential operator and $G: C_0^\infty(E)\to C^\infty(F)$ a Green operator.
\end{definition}

\begin{remark}
 The operators $A$ and $T$ come from the boundary value problems. The
 potential operator $K$ is a solution operator, and the Green operator
 $G$ had to be introduced to obtain an algebra.

Every BdM operator has an order $\mu\in[-\infty,\infty)$ and a type $d\in\N_0$. The order is a generalization of the order of a (pseudo)differential operator, the type is determined by the trace and Green operator and says ``how much restriction to the boundary'' is involved.

BdM operators are locally defined (except the smoothing operators, see
\ref{def_sm}): $\BdM{P}$ is BdM, if and only if for all cutoff functions $\varphi$ and $\psi$ ($\psi=1$ on $\supp\varphi$) with support in a chart the Euclidian operator $\varphi\BdM{P}\psi$ is BdM and if $\varphi\BdM{P}(1-\psi)$ is a BdM operator of order $-\infty$ and type zero.
\end{remark}

For the rest of the section, adopt the situation of Definition \ref{def_BdM}.

\begin{proposition}
Suppose $M$ is compact.\\
Let $\BdM{P}:C^\infty(E)\oplus C^\infty(X)\to C^\infty(F)\oplus C^\infty(Y)$ and $\BdM{Q}:C^\infty(F)\oplus C^\infty(Y)\to C^\infty(G)\oplus C^\infty(Z)$ be BdM operators of order $\mu$ and type $d$ and $\mu'$, $d'$ respectively. Then the composition $\BdM{QP}$ is a BdM operator of order $\mu+\mu'$ and type $\max\{d',d+\mu'\}$.

The BdM operator $\BdM{P}$ of order $\mu\ge 0$ and type $d\le\mu$ is
elliptic if and only if there exists a BdM operator $\BdM{Q}:C^\infty(F)\oplus C^\infty(Y)\to C^\infty(E)\oplus C^\infty(X)$ of order $-\mu$ and type zero so that
\[ \BdM{S}_0:= \BdM{QP}-\1\qquad\text{and}\qquad \BdM{S}_1:=\BdM{PQ}-\1 \]
are of order $-\infty$ and $\BdM{S}_0$ is of type $\mu$, $\BdM{S}_1$ of type zero. $\BdM{Q}$ is called a \em{parametrix} of $\BdM{P}$. Two parametrices differ by an operator of order $-\infty$.

Every differential boundary value problem $\bvp{P}=(A, T):C_0^\infty(E)\to C_0^\infty(F)\oplus C_0^\infty(Y)$ is a Boutet de Monvel operator. If it is elliptic in the Lopatinsky-Shapiro sense, it is also elliptic in the sense of the BdM algebra.
\end{proposition}
\begin{proof}
Compare Schrohe/Schulze \cite{Schrohe-Schulze(1994)} and Rempel/Schulze \cite{Rempel-Schulze(1982)}.
\end{proof}

In the study of elliptic boundary value problems,  the BdM operators of order $-\infty$ are important. These are operators with smooth integral kernels:
\begin{definition}\label{def_sm}
We call $\BdM{P}$ a {\em smoothing Boutet de Monvel operator} (an operator of order $-\infty$ and type $d\ge 0$), if there exist smooth integral kernels \[\begin{split}
 a, g_0 & \in C^\infty(\Hom(p_2^*E,p_1^*F)\bundleover M\times M),\\
g_1,\dots,g_d & \in C^\infty(\Hom(p_2^*E|_{\boundary M}, p_1^*F)\bundleover\, M\times\boundary M),\\
k & \in C^\infty(\Hom(p_2^*X,p_1^*F)\bundleover\, M\times\boundary M),\\
 p & \in C^\infty(\Hom(p_2^*X,p_1^*Y)\bundleover\, \boundary M\times \boundary M),\\
t_1,\dots,t_d & \in C^\infty(\Hom(p_2^* E|_{\boundary M}, p_1^*Y)\bundleover \boundary M\times\boundary M),\\
t_0 & \in C^\infty(\Hom(p_2^*E,p_1^*Y)\bundleover \boundary M\times M),
\end{split}\]
so that $A$ and $p$ have the corresponding integral kernel and for
$F\in C_0^\infty(E)$ and $f\in C_0^\infty(X)$ we have in addition
\[\begin{split}
GF(x) &=\sum_{i=1}^d \int_{\boundary M} g_i(x,y')(\partial_\nu)^{i-1} F(y')dy' +\int_M g_0(x,y)F(y)dy\\
Kf(x) &=\int_{\boundary M}k(x,y')f(y')dy'\\
TF(x') &=\sum_{i=1}^d \int_{\boundary M} t_i(x',y')(\partial_\nu)^{i-1} F(y')dy' + \int_M t_0(x',y)F(y)dy,
\end{split}\]
where $\partial_\nu$ denotes differentiation in inward normal direction.
\end{definition}

\begin{definition} 
Equip $M$  with a Riemannian metric. An operator
$\BdM{P}: C_0^\infty(E)\oplus C_0^\infty(X) \to C^\infty(F)\oplus
C^\infty(Y)$ is called $\epsilon$-local ($\epsilon >0$), if
\[ \supp(\BdM{P}f)\subset \{x\in M;\;d(x,\supp f)<\epsilon\}\quad\forall f\in C_0^\infty .\]
\end{definition}

\begin{proposition}\label{local_ex}
Suppose $M$ is a compact Riemannian manifold and $\epsilon>0$ is given.
Every BdM operator $\BdM{P}$ is the sum of an $\epsilon$-local operator and a smoothing BdM operator of type zero.
\end{proposition}
\begin{proof}
Choose a finite covering of $M$ by balls $\{U_i\}$ of radius $\epsilon/2$. Let $\{\varphi_i\}$ be a subordinate partition of unity and $\psi_i$ cutoff functions with $\psi_i=1$ on $\supp\varphi_i$ and $\supp\psi_i\subset U_i$.
Set
\[\BdM{P}_1:=\sum_i\varphi_i\BdM{P}\psi_i,\quad \BdM{P}_2:=\BdM{P}-\BdM{P}_1 = \sum_i\varphi_i\BdM{P}(1-\psi_i).\]
Then $\BdM{P}_2$ is a smoothing BdM operator of type zero and $\BdM{P}_1$ is $\epsilon$-local. $\BdM{P}_1$ is a BdM operator of same order and type as $\BdM{P}$ \cite[2.3.3.2. Theorem 1]{Rempel-Schulze(1982)}.
\end{proof}

\begin{proposition}
Suppose $M$ is compact. 
Let $\BdM{P}$ be a BdM operator of order $\mu$ and type $d$. If $s>d-1/2$, then $\BdM{P}$ extends to a continuous operator
\[ \BdM{P}: H^s(E)\oplus H^s(X) \to H^{s-\mu}(F)\oplus H^{s-\mu}(Y).\]

\end{proposition}
\begin{proof}
Compare Schrohe/Schulze \cite[2.2.19]{Schrohe-Schulze(1994)}.
\end{proof}

\begin{proposition}
Let $\tilde M\bundleover M$ be a normal Riemannian covering of Riemannian manifolds with covering group $\Gamma$, where $M$ is compact. Suppose the covering is trivial over balls of radius $2\epsilon$. Suppose 
\[\BdM{P}: C^\infty(E)\oplus C^\infty(X)\to C^\infty(F)\oplus C^\infty(Y)\]
is an $\epsilon$-local operator which extends to a bounded operator 
\[\BdM{P}: H^s(E)\oplus H^s(X)\to H^{s-\mu}(F)\oplus H^{s-\mu}(Y).\]
Then $\BdM{P}$ lifts to an operator 
\[\tilde\BdM{P}: C^\infty(\tilde E)\oplus C^\infty(\tilde X)\to C^\infty(\tilde F)\oplus C^\infty(\tilde Y),\]
 which has a bounded extension
\[ \tilde\BdM{P}: H^s(\tilde E)\oplus H^s(\tilde X)\to
H^{s-\mu}(\tilde F)\oplus H^{s-\mu}(\tilde Y) . \]
\end{proposition}
\begin{proof}
Let $\{U_i\}_{i=1,\dots,N}$ be a covering of $M$ by balls of radius
$\epsilon$, let $V_i$ be the corresponding balls of radius
$2\epsilon$. Let $\varphi_i$ be a subordinate covering of unity. This
induces a $\Gamma$-invariant covering
$\{U_{i,\gamma}\}_{\gamma\in\Gamma}$ of $\tilde M$ with subordinate
$\Gamma$-invariant partition of unity $\varphi_{i,\gamma}$. It is
clear how to lift $\BdM{P}$. To check boundedness, let $\mathcal{F}=(F,f)\in C_0^\infty(\tilde E)\oplus C_0^\infty(\tilde X)$ be given. Then (use $\abs{a+b}^2\le 3(\abs{a}^2+\abs{b}^2)$)
\[ \begin{split}
\abs{\tilde\BdM{P}\mathcal{F}}^2_{H^{s-\mu}} &=
\abs{\tilde\BdM{P}\sum_{i,\gamma}\varphi_{i,\gamma}\mathcal
  F}^2_{H^{s-\mu}}\stackrel{3^N}{\le}
\sum_{i=1}^N\abs{\tilde\BdM{P}\sum_{\gamma}\varphi_{i,\gamma}\mathcal{F}}^2_{H^{s-\mu}}\\
&\stackrel{(*)}{=}
\sum_{i,\gamma}\abs{\tilde\BdM{P}\varphi_{i,\gamma}\mathcal{F}}^2_{H^{s-\mu}}
\stackrel{\norm{\BdM{P}}^2}{\le}
\sum_{i,\gamma}\abs{\varphi_{i,\gamma}\mathcal{F}}^2_{H^s} 
\stackrel{\text{Def}}{=}\abs{\mathcal{F}}^2_{H^s} .
\end{split}\]
 $(*)$ holds since $\supp(\varphi_{i,\gamma}) \cap \supp(\varphi_{i,\gamma'})=\emptyset$ if $\gamma\ne\gamma'$.
\end{proof}

Next we compute the trace of sufficiently regularizing BdM
operators. Most important is the fact that the $\Gamma$-trace of a
lift equals the trace of the operator on the base.

\begin{theorem}\label{loc_tr}
Let $\BdM{P}:C^\infty(E)\oplus C^\infty(X)\to C^\infty(E)\oplus
C^\infty(X)$ be a BdM operator of order $-\mu<-4m$ and type $d$
($m=\dim M$). For $s>d-1/2$, $\BdM{P}$ extends to a bounded trace class operator
\[ \BdM{P}: H^s(E)\oplus H^s(X) \to H^s(E)\oplus H^s(X). \]
The value of the trace is independent of $s$.

If $\BdM{P}$ is $\epsilon$-local then its lift $\tilde\BdM{P}:H^s(\tilde E)\oplus H^s(\tilde X)\to H^s(\tilde E)\oplus H^s(\tilde X)$ (defined for $s>d-1/2)$ is $\Gamma$-tr and
\[ \Sp_\Gamma(\tilde\BdM{P}) =\Sp(\BdM{P}). \]
If $-\mu=-\infty$ and $\BdM{P}$ is has  kernels as in  Definition
\ref{def_sm} then explicitly
\[\begin{split}
 \Sp(\BdM{P}) = & \int_M \Sp_{E_x} a(x,x) dx +\int_{\boundary M}\Sp_{X_x'}p(x',x') dx'\\
  & +\int_M \Sp_{E_x} g_0(x,x)dx +\sum_{i=1}^d \int_{\boundary M}\Sp_{E_{x'}}\partial_{\nu,x}^{i-1}p_i(x,y)|_{x=x'=y} dx'
  \end{split}\]
($\Sp_F$ denotes the trace on the finite dimensional vector space $F$; $\partial_\nu$ is differentiation in normal direction).
\end{theorem}

\begin{proof}
The inclusion $H^{s+\mu}\hookrightarrow H^s$ is of trace class by
Theorem \ref{L2_Rellich}). Therefore
$\BdM{P}:H^s\stackrel{\BdM{P}}{\to} H^{s+\mu}\hookrightarrow H^s$ is
of trace class, being the composition of a bounded operator and a trace class operator. If $\mu-4m>s'-s>0$ then
\[\begin{split} 
\Sp(\BdM{P}:H^{s'}\to H^{s'}) = \Sp( H^{s'}\hookrightarrow H^s\stackrel{\BdM{P}}{\to}H^{s+\mu}\hookrightarrow H^{s'})\\
 = \Sp(H^s\stackrel{\BdM{P}}{\to}H^{s+\mu}\hookrightarrow H^{s'}\hookrightarrow{H^s}) = \Sp(\BdM{P}:H^s\to H^s) .
\end{split}\]
Here we used the trace property, noting that $H^{s+\mu}\hookrightarrow H^s$ is trace class. Inductively, the trace is independent of $s$ for arbitrary $s$.

Identical arguments apply to the lift $\tilde\BdM{P}$, replacing trace by $\Gamma$-trace and using Theorem \ref{L2_Rellich}.

Now we come to the explicit computation, and
$\mu=-\infty$. Observe (with the notion of \ref{def_BdM})
$\Sp(\BdM{P})=\Sp(A)+\Sp(G)+\Sp(p)$. Note that $A$ and $p$ are
actually defined on $L^2$. The above argument applies to show that
$\Sp(A:H^s\to H^s)=\Sp(A:L^2\to L^2)$. $A$ is an integral
operator with a smooth kernel and therefore with trace (on $L^2$)
\[ \Sp(A)=\int_M \Sp_{E_x} a(x,x) dx. \]
Similarly $\Sp(p)=\int_{\boundary M}\Sp_{X_{x'}}p(x',x') dx'$.
For the obvious splitting $G=G_0+G_1+\dots+ G_d$, note that each summand is trace class. $G_0$ behaves exactly as $A$ does. For $i>0$, the operator $G_i$ is a composition
\[
H^s(E)\stackrel{\partial_\nu^{i-1}}{\to}H^{s-i+1}(E)\stackrel{\res}{\to}
H^{s-i+1/2}(E|_{\boundary
M})\stackrel{K_i}{\to}H^{\infty}(E)\stackrel{i}{\hookrightarrow}
H^s(E) . \]
Each of the operators is bounded and the inclusion is trace class ($\res$ denotes the restriction to the boundary and $K_i$ is the obvious integral operator with smooth kernel from $E|_{\boundary M}\to E$). Using the trace property and the fact that inclusions of Sobolev spaces commute with differentiation and restriction to the boundary, we see
\[ \Sp(G_i)= \Sp(i\circ \underbrace{\res\circ \partial_\nu^{i-1}\circ
K_i}_{=: P_i}) . \]
Now $P_i$ is an integral operator with smooth kernel on $\boundary M$, namely
\[ P_if(x')=\int_{\boundary M} (\partial^{i-1}_{\nu,x}g_i) (x',y') f(y') dy'. \]
Therefore it extends to a trace class operator on $L^2(E|_{\boundary M})$ with
\[ \Sp(G_i)=\Sp(P_i)=\int_{\boundary M} \Sp_{E_{x'}}(\partial_{\nu,x}^{i-1}g_i(x,y))|_{x=x'=y} dx' .\]
This establishes the formula for $\Sp(\BdM{P})$.

Identical arguments apply to the lift $\tilde\BdM{P}$ as far as follows: \[\Sp_\Gamma(\tilde\BdM{P})=\Sp_\Gamma(\tilde A) + \Sp_\Gamma(\tilde p) + \Sp_\Gamma(\tilde G_0) +\sum_{i=1}^d\Sp_\Gamma(\tilde P_i) ,\]
where each summand is the lift of an integral operator with smooth kernel on
$L^2(E)$, $L^2(X)$ and $L^2(E|_{\boundary M})$, respectively.

Therefore, it remains to show that for this type of operator the $\Gamma$-trace of the lift coincides with the trace on the base: choose a fundamental domain $U\subset \tilde M$ for the covering. Then $L^2(\tilde E)\to l^2(\Gamma)\tensor L^2(E|U)$ is an isometric $\Neumann{N}(\Gamma)$-isomorphism and $\tilde A=\1\tensor A$ (as a lift). By Theorem \ref{tr_prop} (\ref{HS_prod}),
\[ \Sp_\Gamma(\tilde A)=\tr_\Gamma(\1)\cdot\Sp(A)=\Sp(A), \]
and similarly for the other operators.
\end{proof}


\section{Proof of the $L^2$-index theorem}\label{sec_proof}

\begin{situation} Let $\tilde M\bundleover M$ be a normal covering of a compact
manifold with boundary with deck transformation group $\Gamma$. Let $\BdM{P}=(A,T):C_0^\infty(E)\to C_0^\infty(F)\oplus C_0^\infty(Y)$ be an elliptic differential boundary value problem on $M$. Denote its lift to $\tilde M$ with $\tilde\BdM{P}: C_0^\infty(\tilde E)\to C_0^\infty(\tilde F)\oplus C_0^\infty(\tilde Y)$. Suppose $\BdM{P}$ has order $\mu\ge 0$ and type $d\le\mu$.

We have the extension $\BdM{P}:H^\mu(E)\to L^2(F)\oplus L^2(Y)$. 

Let $H_0:L^2(E)\to\ker(\BdM{P})$ be the orthogonal projection onto the kernel, $H_1:L^2(F)\oplus L^2(Y)\to \im(\BdM{P})^\perp$ the orthogonal projection onto the cokernel of $\BdM{P}$. Similarly, let $\bar H_0$ and $\bar H_1$ be the projections onto kernel and cokernel of $\tilde\BdM{P}$.
\end{situation}

We want to prove the $L^2$-index Theorem \ref{ind_theorem} for $\boundary$-manifolds:
\begin{theorem} 
$\dim_\Gamma\ker\tilde\BdM{P}=\Sp_\Gamma(\bar H_0)$ and $\dim_\Gamma\coker\tilde\BdM{P}=\Sp_\Gamma(\bar H_1)$ are finite, and 
\begin{equation*} 
\ind_\Gamma(\tilde\BdM{P}):=\Sp_\Gamma(\bar H_0)-\Sp_\Gamma(\bar H_1)\quad=\quad \ind(\BdM{P})=\Sp(H_0)-\Sp(H_1).
\end{equation*}
\end{theorem}

The idea of the proof is the following: $H_i$ and $\bar H_i$ have in
general nothing to do with each
other. But suppose we could find a bounded liftable "inverse"
$\BdM{Q}$ to $\BdM{P}$. Then the equations
\[ \BdM{P}\BdM{Q}= \1-H_1\quad\text{and}\quad \BdM{Q}\BdM{P}=\1- H_0 \]
could be lifted and we could compare the trace of $H_i$ and $\bar H_i$
directly. This is not possible. We use a parametrix instead:

Let $\BdM{Q}$ be an $\epsilon$-local parametrix of $\BdM{P}$ (use
Proposition \ref{local_ex}) so that
\begin{equation}\label{inv_eq}
  \begin{split}    
    & \BdM{P}\BdM{Q}=\1-\BdM{S}_1,\quad
    \BdM{Q}\BdM{P}=\1-\BdM{S}_0\\
    \implies \quad &
    \tilde\BdM{P}\tilde\BdM{Q}=\1-\tilde\BdM{S}_1,\quad
    \tilde\BdM{Q}\tilde\BdM{P}=\1-\tilde\BdM{S}_0.
\end{split}
\end{equation}
Automatically, $\BdM{S}_0=\1-\BdM{QP}$ and $\BdM{S}_1$ are $\epsilon$-local since the right hand side is.
Note that $\BdM{S}_0$ and $\tilde\BdM{S}_0$ are operators of order $-\infty$ and type $\mu$, whereas $\BdM{S}_1$ and $\tilde \BdM{S}_1$ have order $-\infty$ and type zero.

We know already that $\Sp_\Gamma\tilde\BdM{S}_i=\Sp\BdM{S}_i$ (Theorem \ref{loc_tr}). It remains to show that we can compute the index also in terms of the $\BdM{S}_i$, namely
\begin{equation}\label{s_ind}
 \Sp\BdM{S}_0-\Sp\BdM{S}_1 =\Sp H_0-\Sp H_1
\end{equation}
(and similarly on $\tilde M$). This will be achieved using Theorem \ref{equaltr}. We start with

\begin{proposition}
The image of the projection $H_0:L^2(E)\to L^2(E)$ (i.e.\ the kernel of $\BdM{P}$) is contained in $H^\infty(E)$ and $H_0$ restricts to a bounded operator $H_0:H^s(E)\to H^{s+t}(E)$ for arbitrary $s,t\ge 0$. Especially $H_0:H^s\to H^s$ is trace class for every $s\ge 0$ and the trace is independent of $s$.

The same holds for $\bar H_0$ if we replace $\tr$ by $\tr_\Gamma$.
\end{proposition}
\begin{proof} 
Elliptic regularity and the corresponding a priori estimates (the
theory works as in the compact case, compare \forget{\cite[4.15]{\citediss}
or }\cite[4.14]{Schick(1998d)} for a generalization) imply
that the kernels of $\BdM{P}$ and $\tilde\BdM{P}$ are contained in
every Sobolev space $H^s(E)$ and $H^s(\tilde E)$ respectively, and that the
Sobolev norms on this subspace are equivalent to the $L^2$-norm. This
implies everything if we consider $H_0$ as composition of the bounded
operator $H_0:H^s\to H^{s+4m}$ with the trace class operator
$i:H^{s+4m}\hookrightarrow H^s$ (and similarly for $\bar H_0$).
\end{proof}

Now we can prove equation \ref{s_ind}.
The following computations are formulated only for the lifted
operators. They are valid also on the base with the obvious changes.

Multiplying the equations in \eqref{inv_eq} with $\bar H_1$ from the left and with $\bar H_0$ from the right, we get
\begin{equation}\label{H_S_eq}
 \bar H_1=\bar H_1\tilde\BdM{S}_1 \qquad \bar H_0=\tilde\BdM{S}_0\bar H_0 ,
 \end{equation}
where the equation for $\bar H_0$ is valid on $H^\mu$ and the one for $\bar H_1$ is valid on all of $L^2$ . By multiplication of \eqref{inv_eq} with $\tilde\BdM{P}$ we get on $H^\mu$
\[ \tilde\BdM{P}\tilde\BdM{S}_0 = \tilde\BdM{S}_1\tilde\BdM{P}. \]

Following Atiyah \cite{Atiyah(1976)} we now define 
\begin{equation*}
\bar T_i:=(1-\bar
H_i)\tilde\BdM{S}_i(1-\bar H_i) \quad(i=0,1) .
\end{equation*}
Because of Theorem
\ref{tr_prop} \eqref{tr_comp}
$\bar T_0$ is a $\Gamma$-tr operator on the Hilbert $\Neumann{N}(\Gamma)$-module  $H^\mu$ and $\bar T_1$ is a $\Gamma$-tr operator on the Hilbert $\Neumann{N}(\Gamma)$-module $L^2$. Since $\bar H_i$ are projectors
\[\begin{split} \Sp_\Gamma \bar T_0 &=
\Sp_\Gamma(\tilde\BdM{S}_0(1-\bar
H_0))=\Sp_\Gamma\tilde\BdM{S}_0-\Sp_\Gamma \bar H_0\qquad\text{(use
\eqref{H_S_eq}}) ,\\
\Sp_\Gamma \bar T_1 &= \Sp_\Gamma((1-\bar
H_1)\tilde\BdM{S}_1)=\Sp_\Gamma\tilde\BdM{S}_1-\Sp_\Gamma \bar
H_1\qquad\text{(use \eqref{H_S_eq})} .
\end{split}\]
Therefore,
\[\Sp_\Gamma \bar T_0=\Sp_\Gamma \bar T_1\quad\iff\quad \Sp_\Gamma\tilde\BdM{S}_0-\Sp_\Gamma\tilde\BdM{S}_1=\Sp_\Gamma\bar H_0-\Sp_\Gamma\bar H_1.\]

Next observe
\[ \ker\tilde\BdM{P}\subset \ker \bar T_0;\qquad
\ker\tilde\BdM{P}^*\subset\ker \bar T_1^* ;\]
\[ 
\tilde\BdM{P}\bar T_0  = \tilde\BdM{P}\tilde\BdM{S}_0-\underbrace{\tilde\BdM{P}\bar H_0}_{=0}\tilde\BdM{S}_0-\tilde\BdM{P}\underbrace{\tilde\BdM{S}_0\bar H_0}_{=\bar H_0}+\underbrace{\tilde\BdM{P}\bar H_0}_{=0}\tilde\BdM{S}_0\bar H_0=\tilde\BdM{S}_1\tilde\BdM{P}=\dots=
\bar T_1\tilde\BdM{P} . \]
Application of Proposition \ref{equaltr} yields $\Sp_\Gamma\bar T_0=\Sp_\Gamma\bar T_1$, i.e.\ $\ind_\Gamma\tilde \BdM{P}=\Sp_\Gamma\tilde\BdM{S}_0-\Sp_\Gamma\tilde\BdM{S}_1$. Similarly, $\ind\BdM{P}=\Sp\BdM{S}_0-\Sp\BdM{S}_1$. Now Theorem \ref{loc_tr} appplied to the $\epsilon$-local smoothing operators $\BdM{S}_0$, $\BdM{S}_1$ finishes the proof of Theorem \ref{ind_theorem}.


\section{Index and adjoint boundary value problems}

The purpose of this section is to simplify the index
formula by replacing the cokernel with the kernel of the adjoint.

\begin{theorem}\label{cokernel}
Let $E,F\bundleover M$, $X,Y\bundleover \boundary M$ be Riemannian
vector bundles, $\bvp{P}:=(A,p): C_0^\infty(E)\to C_0^\infty(F)\oplus
C_0^\infty(Y)$ an elliptic differential boundary value problem. If
$\bvp{Q}:=(B,q): C_0^\infty(F)\to C_0^\infty(E)\oplus C_0^\infty(X)$
is adjoint to $(A,p)$ with respect to the Greenian formula
\begin{equation}\label{Green1}
 (Ae,f)_{L^2(F)}- (e,Bf)_{L^2(E)} = (pe,sf)_{L^2(Y)}-(te,qf)_{L^2(X)}
\end{equation}
($e\in C_0^\infty(E), f\in C_0^\infty(F)$ and $t,s$ are auxiliary
boundary differential operators), then 
\[ L^2(F)\oplus L^2(Y)\supset\qquad\im(\bvp{P})^\perp \stackrel{p_1}{\to} L^2(F): (f,y)\mapsto f \]
is an isomorphism onto $\ker(\bvp{Q})$ with inverse
\[ \alpha:\ker(\bvp{Q})\to \im(\bvp{P})^\perp: f\mapsto (f,-sf) .\]
\end{theorem}
\begin{proof}
\cite[13.1]{\citediss}
\end{proof}

Being in the situation of the $L^2$-index Theorem \ref{ind_theorem},
the isomorphism of Theorem \ref{cokernel} is equivariant under the
group operation and $\coker(\tilde\bvp{P})$ is $\Gamma$-isomorphic to $\ker(\tilde\bvp{Q})$. Therefore the index theorem can be stated as follows:
\begin{theorem}\label{ind_theorem+}
Suppose $M$ is a compact boundary manifold with normal covering $\tilde M$ and covering group $\Gamma$. Let $\bvp{P}:=(A,T)$ be an elliptic differential boundary value problem on $M$ with lift $\tilde\bvp{P}$. Let $\bvp{Q}:=(B,S)$ be an adjoint with lift $\tilde\bvp{Q}$. Then
\[ \ind(\bvp{P})= \quad\ind_\Gamma(\tilde\bvp{P})=
\dim_\Gamma(\ker\tilde\bvp{P}) -\dim_\Gamma(\ker\tilde\bvp{Q}) . \]
\end{theorem}

We apply this to compute the Euler characteristic of a
$\boundary$-manifold. Lott/L\"uck \cite{Lott-Lueck(1995)} get the
same result with other methods.

\begin{theorem}\label{Euler}
Suppose $M$ is a compact manifold with boundary $\boundary
M=M_1\amalg M_2$. Let $\tilde M$ be a normal covering of $M$ with
covering group $\Gamma$. Then
\[ \chi(M,M_1) = \sum_p (-1)^p\dim_\Gamma \mathcal{H}^p_{(2)}(\tilde M,\tilde M_1) \]
with $\mathcal{H}^p_{(2)}(\tilde M,\tilde M_1)=\{\omega\in C^\infty(\Lambda^pT\tilde M);\;\abs{\omega}_{L^2}<\infty, d\omega=0=\delta\omega, b_1^*(\omega)=0=b_2^*(*\omega)\}$. ($b_i:\tilde M_i\hookrightarrow \tilde M$ are the inclusions).
\end{theorem}
\begin{proof} To keep notation simple suppose $M_1=\emptyset$.
We known  $\chi(M)=\ind(\bvp{P}^{ev})$, where $\bvp{P}^{ev/odd}$ are
the boundary value problems
\[ (d+\delta, b_2^*\circ *):C^\infty(\Lambda^{ev/odd}TM)\to C^\infty(\Lambda^{odd/ev}TM)\oplus  C^\infty(\Lambda^*T\boundary M). \]
We have the following Greenian formula
  \begin{multline*}
    ((d+\delta)\omega,\eta)_{L^2(M)}=\\
    =(\omega,(\delta+d)\eta)_{L^2(M)}
    \pm \int_{\boundary M} b^*\omega\wedge b^*(*\eta)
    \pm\int_{\boundary M} b^*(\eta)\wedge b^*(*\omega).
\end{multline*}
Theorems \ref{cokernel} and \ref{ind_theorem+}  yield then
\[ \chi(M)=\ind(\bvp{\tilde P}^{ev})=\dim\ker(\bvp{\tilde
  P}^{ev})-\dim\ker(\bvp{\tilde P}^{odd}) . \]
In view of elliptic regularity this is just the claim.
\end{proof}


\bibliographystyle{lueck}
\small\bibliography{literatur}

\begin{thebibliography}{10}
\newcommand{\enquote}[1]{``#1''}

\bibitem{Atiyah(1976)}
{\bf Atiyah, M.}: \enquote{{\em Elliptic operators, discrete groups and von
  {N}eumann algebras\/}}, Ast\'erisque 32, 43--72 (1976)

\bibitem{Atiyah-Bott(1964)}
{\bf Atiyah, M. and Bott, R.}: \enquote{{\em The index problem for manifolds
  with boundary\/}}, in: {\em Differential analysis (papers presented at the
  {B}ombay colloquium 1964)\/},  175--186, Oxford University Press (1964)

\bibitem{Atiyah-Patodi-Singer(1975a)}
{\bf Atiyah, M., Patodi, V.K., and Singer, I.M.}: \enquote{{\em Spectral
  asymmetry and {R}iemannian geometry I\/}}, Math. Proc. Cam. Phil. Soc. 77,
  43--69 (1975)

\bibitem{Dixmier(1969)}
{\bf Dixmier, J.}: \enquote{{\em Les alg\`ebre d'op\'erateurs dans l'espace
  Hilbertien (alg\`ebres de von {N}eumann\/}}, Gauthier-Villars (1969)

\bibitem{Dodziuk(1979)}
{\bf Dodziuk, J.}: \enquote{{\em $L^2$-harmonic forms on rotationally symmetric
  {R}iemannian manifolds\/}}, Proc. of the AMS 77, 395--400 (1979)

\bibitem{Donnelly-Xavier(1984)}
{\bf Donnelly, H. and Xavier, F.}: \enquote{{\em On the differential form
  spectrum of negatively curved {R}iemannian manifolds\/}}, Amer. J. of Math.
  106, 169--185 (1984)

\bibitem{Lott-Lueck(1995)}
{\bf Lott, J. and L{\"u}ck, W.}: \enquote{{\em $L^2$-topological invariants of
  $3$-manifolds\/}}, Inventiones Mathematicae 120, 15--60 (1995)

\bibitem{Palais(1965)}
{\bf Palais, R.S.}: \enquote{{\em Seminar on the Atiyah-Singer index
  theorem\/}}, vol.~57 of {\em Annals of mathematics studies\/}, Princeton
  University Press (1965)

\bibitem{Ramachandran(1993)}
{\bf Ramachandran, M.}: \enquote{{\em Von {N}eumann index theorems for
  manifolds with boundary\/}}, Journal of Differential Geometry 38, 315--349
  (1993)

\bibitem{Rempel-Schulze(1982)}
{\bf Rempel, S. and Schulze}: \enquote{{\em Index theory of elliptic boundary
  value problems\/}}, Akademie Verlag, Berlin (1982)

\bibitem{Schick(1994)}
{\bf Schick, T.}: \enquote{{\em Sobolev Spaces and $L^2$-Hodge-De Rham Theorem
  for Infinite Coverings of Manifolds with Boundary\/}}, Diplomarbeit, Johannes
  {G}utenberg-Universit{\"a}t Mainz (1994)

\bibitem{Schick(1996)url}
{\bf Schick, T.}: \enquote{{\em Analysis on $\partial$-manifolds of bounded
  geometry, {H}odge-de {R}ham isomorphism and $L^2$-index theorem\/}}, Shaker,
  Aachen, (Dissertation, Mainz),
  http://wwwmath.uni-muenster.de/math/inst/reine/\\ inst/lueck/publ/schick/disssc%
hick.html (1996)

\bibitem{Schick(1998d)}
{\bf Schick, T.}: \enquote{{\em Geometry and Analysis on $\partial$-manifolds
  of bounded geometry\/}}, preprint, M{\"u}nster (1998)

\bibitem{Schrohe-Schulze(1994)}
{\bf Schrohe, E. and Schulze, B.-W.}: \enquote{{\em Boundary value problems in
  Boutet de Monvel's algebra for manifolds with conical singularities I\/}},
  in: {\em Pseudodifferential operators and mathematical physics. Advances in
  partial differential equations 1\/},  97--209, Akademie Verlag (1994)

\bibitem{Taylor(1981)}
{\bf Taylor, M.E.}: \enquote{{\em Pseudodifferential operators\/}}, Princeton
  University Press (1981)

\end{thebibliography}

\noindent
{\small Thomas Schick\\
FB Mathematik, Universit\"at M\"unster\\
Einsteinstr.~62, 48149 M\"unster, Germany\\
e-mail: thomas.schick@math.uni-muenster.de}
\end{document}